\documentclass[12pt]{article}
\usepackage{amsmath,amssymb,amsbsy,amsfonts,amsthm,latexsym,amsopn,amstext,
                                     amsopn,amstext,amsxtra,euscript,amscd}
\begin{document}

\newfont{\teneufm}{eufm10}
\newfont{\seveneufm}{eufm7}
\newfont{\fiveeufm}{eufm5}
%
%
\newfam\eufmfam
                     \textfont\eufmfam=\teneufm
\scriptfont\eufmfam=\seveneufm
                     \scriptscriptfont\eufmfam=\fiveeufm

%
%
%


\def\bbbr{{\rm I\!R}} 
\def\bbbc{{\rm I\!C}} 
\def\bbbm{{\rm I\!M}}
\def\bbbn{{\rm I\!N}} 
\def\bbbf{{\rm I\!F}}
\def\bbbh{{\rm I\!H}}
\def\bbbk{{\rm I\!K}}
\def\bbbl{{\rm I\!L}}
\def\bbbp{{\rm I\!P}}
\newcommand{\lcm}{{\rm lcm}}
\def\bbbone{{\mathchoice {\rm 1\mskip-4mu l} {\rm 1\mskip-4mu l}
{\rm 1\mskip-4.5mu l} {\rm 1\mskip-5mu l}}}
\def\bbbc{{\mathchoice {\setbox0=\hbox{$\displaystyle\rm C$}\hbox{\hbox
to0pt{\kern0.4\wd0\vrule height0.9\ht0\hss}\box0}}
{\setbox0=\hbox{$\textstyle\rm C$}\hbox{\hbox
to0pt{\kern0.4\wd0\vrule height0.9\ht0\hss}\box0}}
{\setbox0=\hbox{$\scriptstyle\rm C$}\hbox{\hbox
to0pt{\kern0.4\wd0\vrule height0.9\ht0\hss}\box0}}
{\setbox0=\hbox{$\scriptscriptstyle\rm C$}\hbox{\hbox
to0pt{\kern0.4\wd0\vrule height0.9\ht0\hss}\box0}}}}
\def\bbbq{{\mathchoice {\setbox0=\hbox{$\displaystyle\rm
Q$}\hbox{\raise
0.15\ht0\hbox to0pt{\kern0.4\wd0\vrule height0.8\ht0\hss}\box0}}
{\setbox0=\hbox{$\textstyle\rm Q$}\hbox{\raise
0.15\ht0\hbox to0pt{\kern0.4\wd0\vrule height0.8\ht0\hss}\box0}}
{\setbox0=\hbox{$\scriptstyle\rm Q$}\hbox{\raise
0.15\ht0\hbox to0pt{\kern0.4\wd0\vrule height0.7\ht0\hss}\box0}}
{\setbox0=\hbox{$\scriptscriptstyle\rm Q$}\hbox{\raise
0.15\ht0\hbox to0pt{\kern0.4\wd0\vrule height0.7\ht0\hss}\box0}}}}
\def\bbbt{{\mathchoice {\setbox0=\hbox{$\displaystyle\rm
T$}\hbox{\hbox to0pt{\kern0.3\wd0\vrule height0.9\ht0\hss}\box0}}
{\setbox0=\hbox{$\textstyle\rm T$}\hbox{\hbox
to0pt{\kern0.3\wd0\vrule height0.9\ht0\hss}\box0}}
{\setbox0=\hbox{$\scriptstyle\rm T$}\hbox{\hbox
to0pt{\kern0.3\wd0\vrule height0.9\ht0\hss}\box0}}
{\setbox0=\hbox{$\scriptscriptstyle\rm T$}\hbox{\hbox
to0pt{\kern0.3\wd0\vrule height0.9\ht0\hss}\box0}}}}
\def\bbbs{{\mathchoice
{\setbox0=\hbox{$\displaystyle     \rm S$}\hbox{\raise0.5\ht0\hbox
to0pt{\kern0.35\wd0\vrule height0.45\ht0\hss}\hbox
to0pt{\kern0.55\wd0\vrule height0.5\ht0\hss}\box0}}
{\setbox0=\hbox{$\textstyle        \rm S$}\hbox{\raise0.5\ht0\hbox
to0pt{\kern0.35\wd0\vrule height0.45\ht0\hss}\hbox
to0pt{\kern0.55\wd0\vrule height0.5\ht0\hss}\box0}}
{\setbox0=\hbox{$\scriptstyle      \rm S$}\hbox{\raise0.5\ht0\hbox
to0pt{\kern0.35\wd0\vrule height0.45\ht0\hss}\raise0.05\ht0\hbox
to0pt{\kern0.5\wd0\vrule height0.45\ht0\hss}\box0}}
{\setbox0=\hbox{$\scriptscriptstyle\rm S$}\hbox{\raise0.5\ht0\hbox
to0pt{\kern0.4\wd0\vrule height0.45\ht0\hss}\raise0.05\ht0\hbox
to0pt{\kern0.55\wd0\vrule height0.45\ht0\hss}\box0}}}}
\def\bbbz{{\mathchoice {\hbox{$\sf\textstyle Z\kern-0.4em Z$}}
{\hbox{$\sf\textstyle Z\kern-0.4em Z$}}
{\hbox{$\sf\scriptstyle Z\kern-0.3em Z$}}
{\hbox{$\sf\scriptscriptstyle Z\kern-0.2em Z$}}}}
\def\ts{\thinspace}

\newtheorem{theorem}{Theorem}
\newtheorem{lemma}[theorem]{Lemma}
\newtheorem{claim}[theorem]{Claim}
\newtheorem{cor}[theorem]{Corollary}
\newtheorem{prop}[theorem]{Proposition}
\newtheorem{definition}{Definition}
\newtheorem{question}[theorem]{Open Question}

\def\squareforqed{\hbox{\rlap{$\sqcap$}$\sqcup$}}
\def\qed{\ifmmode\squareforqed\else{\unskip\nobreak\hfil
\penalty50\hskip1em\null\nobreak\hfil\squareforqed
\parfillskip=0pt\finalhyphendemerits=0\endgraf}\fi}

\def\cA{{\mathcal A}}
\def\cB{{\mathcal B}}
\def\cC{{\mathcal C}}
\def\cD{{\mathcal D}}
\def\cE{{\mathcal E}}
\def\cF{{\mathcal F}}
\def\cG{{\mathcal G}}
\def\cH{{\mathcal H}}
\def\cI{{\mathcal I}}
\def\cJ{{\mathcal J}}
\def\cK{{\mathcal K}}
\def\cL{{\mathcal L}}
\def\cM{{\mathcal M}}
\def\cN{{\mathcal N}}
\def\cO{{\mathcal O}}
\def\cP{{\mathcal P}}
\def\cQ{{\mathcal Q}}
\def\cR{{\mathcal R}}
\def\cS{{\mathcal S}}
\def\cT{{\mathcal T}}
\def\cU{{\mathcal U}}
\def\cV{{\mathcal V}}
\def\cW{{\mathcal W}}
\def\cX{{\mathcal X}}
\def\cY{{\mathcal Y}}
\def\cZ{{\mathcal Z}}

\newcommand{\comm}[1]{\marginpar{%
\vskip-\baselineskip 
\raggedright\footnotesize
\itshape\hrule\smallskip#1\par\smallskip\hrule}}





\hyphenation{re-pub-lished}

\def\ord{{\mathrm{ord}}}
\def\Nm{{\mathrm{Nm}}}
\renewcommand{\vec}[1]{\mathbf{#1}}

\def \F{{\bbbf}}
\def \L{{\bbbl}}
\def \K{{\bbbk}}
\def \Z{{\bbbz}}
\def \N{{\bbbn}}
\def \Q{{\bbbq}}
\def\E{{\mathbf E}}
\def\G{{\mathcal G}}
\def\O{{\mathcal O}}
\def\cS{{\mathcal S}}
\def \R{{\bbbr}}
\def\Fp{\F_p}
\def \fp{\Fp^*}
\def\\{\cr}
\def\({\left(}
\def\){\right)}
\def\fl#1{\left\lfloor#1\right\rfloor}
\def\rf#1{\left\lceil#1\right\rceil}

\def\Zm{\Z_m}
\def\Zt{\Z_t}
\def\Zp{\Z_p}
\def\Um{\cU_m}
\def\Ut{\cU_t}
\def\Up{\cU_p}

\def\ep{{\mathbf{e}}_p}

\def \Diam{{\mathrm {Diam}}}
\def \Dist{{\mathrm {Dist}}}

\def\LC{{\cL}_{C,\cF}(Q)}
\def\LCn{{\cL}_{C,\cF}(nG)}
\def\Mrs{\cM_{r,s}\(\F_p\)}

\def\Fbar{\overline{\F}_q}
\def\Fn{\F_{q^n}}
\def\En{\E(\Fn)}

\def\mand{\qquad \mbox{and} \qquad}

\def\vf{\vec{f}}


\title{On Point Sets in Vector Spaces over Finite Fields 
That Determine  Only Acute Angle Triangles}

\author{ 
\sc{Igor E.~Shparlinski}\\
{Department of Computing}\\
{Macquarie University}\\
{Sydney, NSW 2109, Australia}\\
{\tt igor@ics.mq.edu.au}
}

\pagenumbering{arabic}

\date{}
\maketitle

\begin{abstract} For three points $\vec{u}$,
$\vec{v}$ and $\vec{w}$ in the $n$-dimensional
space $\F_q^n$
over the finite field $\F_q$ of $q$ elements  
we give a natural interpretation of an acute angle 
triangle defined by this points. We obtain an 
upper bound on the 
size of a set $\cZ$ such that 
all triples of distinct points 
$\vec{u}, \vec{v}, \vec{w} \in \cZ$ 
define acute angle triangles. A similar question 
in the real space $\cR^n$ dates back to P.~Erd{\H o}s
and has been studied by several authors. 
\end{abstract}


\section{Introduction}

Recent remarkable results of  Bourgain,   Katz and Tao~\cite{BKT}
on the sum-product problem in finite fields  have stimulated 
a series of studies of finite field analogues  of
classical combinatorial and discrete geometry  problems,
see~\cite{BouGar,CHIKR,Dvir,Gar,HaIo,HIKR,HaIoSo,IosKoh,
IosRud,IoShXi,IosSen,KatzShen,Shen,Vinh1,Vinh2,Vinh3,Vinh4} 
and references therein.

Here we extend the scope of such problems and consider 
the question about the largest cardinality of a 
set of   points in the $n$-dimensional space
over a finite field such that every  triples of 
distinct points of this set
defines an acute angle triangle. 
We note that a similar question 
in the Euclidean space $\cR^n$ dates back to P.~Erd{\H o}s
and has been studied by several authors, see~\cite{AckBen-Zw}.

Certainly the notion of an acute triangle (or angle) is 
not immediately obvious in vector spaces over finite 
fields. Here we use the ``rational'' interpretation of
trigonometry invented 
by Wildberger~\cite{Wild} to extend this notion to finite fields. 

To motivate our definition, we note that in the triangle defined by 
three distinct vectors $\vec{u}, \vec{v}, \vec{w} \in \R^n$, 
the vertex at $\vec{u}$ has an acute  angle 
if and  only if
$$
\|\vec{u}-\vec{v}\|^2 + \|\vec{u}-\vec{w}\|^2 -
\|\vec{v}-\vec{w}\|^2 > 0, 
$$
where $\|\vec{x}\|$
is the Euclidean norm of $\vec{x} \in \R^n$.

We now identify positive elements of a finite field $\F_q$ 
of $q$ elements with quadratic residues in $\F_q$ 
and say that 
in the triangle defined by 
three distinct vectors 
$$\vec{u} = (u_1, \ldots, u_n),\ \vec{v}=(v_1, \ldots, v_n), \
\vec{w}=(w_1, \ldots, w_n) \in \F_q^n
$$ 
that the vertex at $\vec{u}$ has an acute  angle 
if any only if
$$
\Delta(\vec{u}, \vec{v}, \vec{w}) = 
\sum_{i=1}^n \((u_i-v_i)^2 + (u_i-w_i)^2-(v_i-w_i)^2\)
$$
is a quadratic residue in $\F_q$.

Since in the field of even characteristic we always have 
$\Delta(\vec{u}, \vec{v}, \vec{w}) = 0$, this definition
makes sense only if $q$ is odd. 

We also remark that 
\begin{equation}
\label{eq:D expand}
\begin{split}
\Delta(\vec{u}, \vec{v}, \vec{w})  &
= 2 \sum_{i=1}^n \(u_i^2 - u_i v_i - u_i w_i + v_iw_i\) = 
\\
& = 2 \(\vec{u} \cdot \vec{u} - \vec{u} \cdot \vec{v} -\vec{u}\cdot \vec{w}+\vec{v}\cdot\vec{w}\)= 2(\vec{u}- \vec{v})\cdot(\vec{u}- \vec{w}),
\end{split} 
\end{equation}
where   $\vec{a}\cdot \vec{b}$ denotes the inner product
of $\vec{a}, \vec{b} \in \F_q^n$.
Thus, 
if $q$ is odd then $\Delta(\vec{u}, \vec{v}, \vec{w}) = 0$ 
if an only if  
$(\vec{u}- \vec{v})\cdot(\vec{u}- \vec{w}) = 0$, which correspond
to the orthogonality at $\vec{u}$   and thus to the  Pythagoras theorem. 

Let $N(n,q)$ be the largest possible  cardinality of a set 
$\cZ \subseteq \F_q$ such that 
all triples of distinct points 
$\vec{u}, \vec{v}, \vec{w} \in \cZ$ 
define acute angle triangles.

We remark  that~\cite[Theorem~1.1]{HaIo}
immediately implies that 
\begin{equation}
\label{eq:N gen}
N(n,q) = O(q^{(n+1)/2}), 
\end{equation}
where the implied constant depends only on $n$. 
In general, we do not know how to improve this
bound. However for $n=2$ we obtain a stronger estimate.

\begin{theorem}
\label{thm: Main}
For a sufficiently large  odd $q$,
$$
N(2,q) \le 2 q^{4/3}. 
$$
\end{theorem}

\section{Additive Character Sums}

Let $\Psi$ be the set of all additive characters of $\F_q$ 
and let $\Psi^*\subset \Psi$ be the set of all nonprincipal 
characters, see~\cite[Section~11.1]{IwKow} for basic properties of additive characters.
In particular, we also recall the identity
\begin{equation}
\label{eq:ident}
 \sum_{\psi \in \Psi} \psi(z) =
\left\{ \begin{array}{ll}
q,& \quad \text{if}\ z=0,\\
0,& \quad \text{otherwise},
\end{array} \right.
\end{equation}
see~\cite[Section~11.1]{IwKow}. 

For an additive character $\psi \in \Psi$ and $\alpha \in \F_q$,  
we define the Gauss sum
$$
G_\psi(\alpha) = \sum_{z \in \F_q} \psi(\alpha z^2) = \sum_{z \in \F_q} \chi(z)\psi(\alpha z)
\chi(\alpha) \sum_{z \in \F_q} \chi(z)\psi(  z) = \chi(\alpha)G_\psi(1), 
$$
where $\chi$ is the quadratic character in $\F_q$ (which exists  since $q$ is odd), 
and recall that 
\begin{equation}
\label{eq:Gauss}
|G_\psi| = q^{1/2},
\end{equation}
for $\psi \in \Psi^*$ and $\alpha \in \F_q^*$, 
see~\cite[Proposition~11.5]{IwKow}.

Finally, given a set  $\cZ \subseteq \F_q^n$, we define the triple character sum
$$
S_\psi(\cZ) = \sum_{\vec{u}, \vec{v}, \vec{w}\in \cZ}
\psi \(\Delta(\vec{u}, \vec{v}, \vec{w})\).
$$
Although we use our result on $S_\psi(\cZ)$ only in the case of 
$n=2$, here we present it in full generality as it may have some 
other applications.

\begin{lemma}
\label{lem: S sum}
For any  $\psi \in \Psi^*$ and a set $\cZ \subseteq \F_q^n$, we have
$$
|S_\psi(\cZ)|^2\le  \# \cZ q^n
 \sum_{\substack{\vec{v}, \vec{w}, \vec{x}, \vec{y}\in \cZ\\ 
 \vec{v} +\vec{w}= \vec{x}+\vec{y}}} 
 \psi\(2\(\vec{v}\cdot\vec{w}-\vec{x}\cdot \vec{y}\)\).
 $$
\end{lemma}

\begin{proof}
We have
$$
|S_\psi(\cZ)| \le   \sum_{\vec{u}\in \cZ}
\left| \sum_{ \vec{v}, \vec{w}\in \cZ}
\psi \(\Delta(\vec{u}, \vec{v}, \vec{w})\)\right|.
$$
Hence, recalling~\eqref{eq:D expand}, we  derive
$$
|S_\psi(\cZ)| \le   \sum_{\vec{u}\in \cZ}
\left| \sum_{ \vec{v}, \vec{w}\in \cZ}
\psi \(-2 \(\vec{u} \cdot \(\vec{v} +\vec{w}\)-\vec{v}\cdot\vec{w}\)\)\right|.
$$ 
Note that since $\psi(-z) = \overline{\psi(z)}$, we can replace $-2$ with $2$.
By the Cauchy inequality
\begin{eqnarray*}
\lefteqn{
|S_\psi(\cZ)|^2 \le  \#\cZ \sum_{\vec{u}\in \cZ}
\left| \sum_{ \vec{v}, \vec{w}\in \cZ}
\psi \(2 \(\vec{u} \cdot \(\vec{v} +\vec{w}\)-\vec{v}\cdot\vec{w}\)\)\right|^2}\\
 & & \qquad  \le  \#\cZ \sum_{\vec{u}\in \F_q^n}
\left| \sum_{ \vec{v}, \vec{w}\in \cZ}
\psi \(2 \(\vec{u} \cdot \(\vec{v} +\vec{w}\)-\vec{v}\cdot\vec{w}\)\)\right|^2\\
 & & \qquad  =  \#\cZ \sum_{\vec{u}\in \F_q^n}
 \sum_{ \vec{v}, \vec{w}, \vec{x}, \vec{y}\in \cZ}
\psi \(2 \(\vec{u} \cdot \(\vec{v} +\vec{w}- \vec{x} - \vec{y}\)
-\vec{v}\cdot\vec{w}+ \vec{x}\cdot \vec{y}\)\)\\
 & & \qquad  =  \#\cZ
 \sum_{ \vec{v}, \vec{w}, \vec{x}, \vec{y}\in \cZ}
 \psi \(2 \( \vec{x}\cdot \vec{y}-\vec{v}\cdot\vec{w}\)\)
  \sum_{\vec{u}\in \F_q^n}
\psi \(2 \vec{u} \cdot \(\vec{v} +\vec{w}- \vec{x} - \vec{y}\)\).
\end{eqnarray*}
Finally, changing the order of summation, and replacing $2 \vec{u}$ 
with $\vec{u}$, we obtain
$$
|S_\psi(\cZ)|^2 \le \cZ
 \sum_{ \vec{v}, \vec{w}, \vec{x}, \vec{y}\in \cZ}
 \psi \(2 \( \vec{x}\cdot \vec{y}-\vec{v}\cdot\vec{w}\)\)
  \sum_{\vec{u}\in \F_q^n}
\psi \(\vec{u} \cdot \(\vec{v} +\vec{w}- \vec{x} - \vec{y}\)\).
$$

By the orthogonality property of additive characters~\eqref{eq:ident},
we see that the inner sum vanishes if and only if  
$$
\vec{v} +\vec{w}- \vec{x} - \vec{y} = 0
$$
in which case it equals $q^n$. Now renaming the variables
$(\vec{v},\vec{w}) \leftrightarrow (\vec{x},\vec{y})$, we conclude the proof. 
\end{proof}

\section{Proof of Theorem~\ref{thm: Main}}

Assume that for set $\cZ \subseteq \F_q^2$ all
triples of distinct vectors $\vec{u}, \vec{v}, \vec{w} \in \cZ$ 
define acute angle triangles. Fix an arbitrary quadratic 
non-residue $\alpha\in \F_q$. 
Then we see that the equation 
$$
\Delta(\vec{u}, \vec{v}, \vec{w}) = \alpha z^2, \qquad
\vec{u}, \vec{v}, \vec{w} \in \cZ, \ z \in \F_q,
$$
has at most 
\begin{equation}
\label{eq:T upper}
T \le \(\# \cZ\)^2 
\end{equation}
solutions 
(which come only from the triples
$\vec{u}, \vec{v}, \vec{w} \in \cZ$ with $\Delta(\vec{u}, \vec{v}, \vec{w})=0$,
that is, when $\vec{u} = \vec{v}$ or $\vec{u} = \vec{w}$).  

On the other hand, from the orthogonality property
of characters~\eqref{eq:ident}, we obtain
$$
T = \sum_{\vec{u}, \vec{v}, \vec{w} \in \cZ} \
\sum_{z \in \F_q} \frac{1}{q}  \sum_{\psi \in \Psi}
\psi\(\Delta(\vec{u}, \vec{v}, \vec{w}) - \alpha z^2\)
= \frac{1}{q}  \sum_{\psi \in \Psi} G_\psi(-\alpha)S_\psi(\cZ).
$$
The term corresponding to the principal character $\psi = \psi_0$
is equal to $\(\# \cZ\)^3$. Thus, recalling~\eqref{eq:Gauss}, we 
obtain 
\begin{equation}
\label{eq:T prelim}
\left|T - \(\# \cZ\)^3\right| \le q^{-1/2}  R ,
\end{equation}
where 
$$
R  = \sum_{\psi \in \Psi^*} |S_\psi(\cZ)|.
$$
Now, by the Cauchy inequality
$$
R^2  = q \sum_{\psi \in \Psi^*} |S_\psi(\cZ)|^2.
$$
Thus using Lemma~\ref{lem: S sum} and then exteding the 
summation to all $\psi \in \Psi$, we deduce
$$
R^2  \le \# \cZ q^3
\sum_{\psi \in \Psi } \sum_{\substack{\vec{v}, \vec{w}, \vec{x}, \vec{y}\in \cZ\\ 
 \vec{v} +\vec{w}= \vec{x}+\vec{y}}} 
 \psi\(2\(\vec{v}\cdot\vec{w}-\vec{x}\cdot \vec{y}\)\). 
$$ 
Changing the order of summation and using~\eqref{eq:ident} again,
we obtain 
\begin{equation}
\label{eq:R and W}
R^2  = \# \cZ q^3
 \sum_{\substack{\vec{v}, \vec{w}, \vec{x}, \vec{y}\in \cZ\\ 
 \vec{v} +\vec{w}= \vec{x}+\vec{y}}} \sum_{\psi \in \Psi }
 \psi\(2\(\vec{v}\cdot\vec{w}-\vec{x}\cdot \vec{y}\)\) =\# \cZ q^4 W,
\end{equation}
where $W$ is the  number of solutions to the system of equations
$$
 \vec{v} +\vec{w}= \vec{x}+\vec{y} \mand \vec{v}\cdot\vec{w}=\vec{x}\cdot \vec{y}
$$
in $\vec{v}, \vec{w}, \vec{x}, \vec{y}\in \cZ$, which is the same
as the  number of solutions to the   equation
$$
\vec{v}\cdot\vec{w} = \vec{x}\cdot\(\vec{v} +\vec{w}- \vec{x}\)
$$
in $\vec{v}, \vec{w}, \vec{x}\in \cZ$. Clearly, when 
$\vec{v}, \vec{w}$ and one component of $\vec{x}$ are fixed,
we obtain a nontrivial quadratic equation over $\F_q$ 
for the other component of $\vec{x}$. Therefore 
$$
W \le 2 \(\# \cZ\)^2 q.
$$
Substituting in~\eqref{eq:R and W}, we obtain 
$$
R^2 \le 2 \(\# \cZ\)^3 q^5.
$$
In turn, inserting this estimate in~\eqref{eq:T prelim} yields
\begin{equation}
\label{eq:T asymp}
\left|T - \(\# \cZ\)^3\right| \le \sqrt{2} \(\# \cZ\)^{3/2}  q^2.
\end{equation}
If $\# \cZ< 2 q^{4/3}$ then there is nothing to proof.
Otherwise,
$$
  \sqrt{2} \(\# \cZ\)^{3/2}  q^2 \le \frac{1}{2}  \(\# \cZ\)^3
$$
thus, by~\eqref{eq:T asymp} we obtain 
$$
T \ge \frac{1}{2}  \(\# \cZ\)^3
$$
which contradicts to~\eqref{eq:T upper},
provided that $q$ is large enough. 
 
\section{Remarks}

Unfortunately the method of this paper, although works for any $n$,
leads to a bound which is that same as~\eqref{eq:N gen} for 
$n =3$ and is even weaker than~\eqref{eq:N gen} for 
$n \ge 4$. 

Furthermore, using the  bound 
$$|S_\psi(\cZ)| \le  \(\# \cZ\)^{2} q^{n/2}
$$
(which is immediate from Lemma~\ref{lem: S sum})
in the argument 
of the proof of Theorem~\ref{thm: Main}
one can recover the bound~\eqref{eq:N gen},
but  it does not seem to
give anything stronger that this.

An alternative way to estimating $N(n,q)$ is via bounds 
of quadratic character sums
$$
T_\chi(\cZ) = \sum_{\vec{u}, \vec{v}, \vec{w}\in \cZ}
\chi \(\Delta(\vec{u}, \vec{v}, \vec{w})\).
$$
Using the same approach (via Cauchy inequality 
and extending summation over $\vec{u}\in \cZ$ to the full 
space $\F_q^n$) as in the proof of Lemma~\ref{lem: S sum}, we obtain 
$$
|T_\chi(\cZ)|^2 = \# \cZ \sum_{ \vec{v}, \vec{w}, \vec{x}, \vec{y}\in \cZ}
\sum_{\vec{u}\in \F_q^n}
\chi \(\Delta(\vec{u}, \vec{v}, \vec{w}) \Delta(\vec{u}, \vec{x}, \vec{y})\).
$$
It is natural to conjecture that the inner sums admits 
a square root estimate and thus is $O(q^{n/2})$ 
unless $(\vec{v}, \vec{w})$ is a permutation of $(\vec{x}, \vec{y})$.
One can derive that $N(n,q) =  O(q^{n/2})$ from such a hypothetical 
bound. Unfortunately the highest form of the polynomial 
$$
F_{ \vec{v}, \vec{w}, \vec{x}, \vec{y}} (\vec{U}) = 
\Delta(\vec{U}, \vec{v}, \vec{w}) \Delta(\vec{U}, \vec{x}, \vec{y})
\in \F[\vec{U}]
$$
is singular, so known analogues of the Deligne bound
for multivariate character sums, see~\cite{Katz,Per},
do not apply. 

We recall that in $\R^n$, the largest number of vectors 
such that each  three of them 
define an   acute  angle  triangle 
is bounded by a function of $n$, see~\cite{AckBen-Zw}.
Although our bounds seem to be much higher that the
true order of magnitude  of $N(n,q)$, we observe that 
$\limsup_{p\to \infty} N(n,p) = \infty$. Indeed, 
by the result of  Graham and   Ringrose~\cite{GrRi},
there is an absolute constant $C > 0$ 
such that for infinitely many primes $p$  all
nonnegative integers $z \le C \log p \log \log \log p$
are quadratic residues modulo $p$. 
Thus, for each such $p$ and an appropriate constant $c$, 
for the set 
$\cZ \subseteq \F_p$ formed by vectors with conponents
in the interval $[1, c\sqrt{\log p \log \log \log p}]$
we have
$$
1 \le \Delta(\vec{u}, \vec{v}, \vec{w}) \le  C \log p \log \log \log p
$$
for any pairwise distinct 
$\vec{u}, \vec{v}, \vec{w} \in \cZ$ and thus $\Delta(\vec{u}, \vec{v}, \vec{w})$
is a quadratic residue. This implies 
$$
\limsup_{p\to \infty} \frac{N(n,p)}{\(\log p \log \log \log p\)^{n/2}} > 0. 
$$

\section*{Acknowledgements}

The authors is grateful to  Norman Wildberger
for his suggestion of an $\F_q$-interpretation 
of acute angle triangles, which has  been used here. 

The authors would also like to thank  Nick Katz for a 
useful discussion 
of some issues related to possible estimates 
of the quadratic character sums $T_\chi(\cZ)$.

During the preparation of this paper, the author was
supported in part by ARC grant DP0556431.

\end{document}